\documentclass[10pt,reqno]{amsart}
\usepackage{amssymb,version,graphicx,fancybox,pifont}
\usepackage{multirow}
\usepackage{url,hyperref}
\usepackage{mathrsfs}
\usepackage{amsmath}
\usepackage{subfigure}
\usepackage{enumerate}
\usepackage{booktabs}
\usepackage{braket}
\usepackage{cleveref}
\usepackage{bm}
\catcode`\@=11
\@addtoreset{equation}{section}   % Makes \section reset 'equation' counter.

\@addtoreset{table}{section}   % Makes \section reset 'table' counter.
\renewcommand\thefigure{\thesection.\@arabic\c@figure}
\@addtoreset{figure}{section}   % Makes \section reset 'figure' counter.
\renewcommand\thetable{\thesection.\@arabic\c@table}

 \newcommand{\new}{\newcommand*}
 \new{\rnew}{\renewcommand*}
  \new{\newe}{\newenvironment*}
 \new{\stl}{\setlength}
 \stl{\arraycolsep}{0.5mm}

\newtheorem{thm}{\bf Theorem}

\newenvironment{theorem}{\begin{thm}} {\end{thm}}
\newtheorem{lmm}{\bf Lemma}

\newenvironment{lemma}{\begin{lmm}}{\end{lmm}}
\theoremstyle{remark}

\theoremstyle{definition}

\usepackage{multirow}
\newcommand {\bgeq}[1]{\begin{equation}\label{#1}}
\newcommand \edeq {\end{equation}}
\newcommand \bgth {\begin{theorem}\label}
\newcommand \edth {\end{theorem}}
\newcommand \bglm {\begin{lemma}\label}
\newcommand \edlm {\end{lemma}}
\newcommand {\bgar}[1]{\begin{array}{#1}}
\newcommand {\edar}{\end{array}}

\usepackage{algorithm, algorithmic}
\title[solving Monge-Amp\`{e}re equation by using Legendre-KAN method]{solving the fully nonlinear Monge-Amp\`{e}re equation using the Legendre-Kolmogorov-Arnold Network method}
\author[B. Hu, L. Jin and Z. Li]{}

\thanks{$^\dag$Corresponding author  email:  zxli@shnu.edu.cn. \\ \indent The work was partially supported by the NNSF of China (No.
12271366, 11871043, 12171322), the NSF of Shanghai, China
(No.21ZR1447200, 22ZR1445500) and the Science and Technology Innovation Plan of Shanghai (No.23JC1403200).}

\begin{document} \maketitle
%% Enter the first author's name and address:
\centerline{Bingcheng Hu$^{1}$, Lixiang Jin$^{2}$, Zhaoxiang Li$^{2,\dag}$}
%\scshape
%\medskip
%{\footnotesize
% %% please put the address of the first author
 \centerline{1. Department of Electronic Engineering, Shanghai Normal University, Shanghai,
200234, P.R. China}
\centerline{2. Department of  Mathematics, Shanghai Normal University, Shanghai,
200234, P.R. China}
 %% Do not forget to end the {\footnotesize by the sign }

%\medskip

%\centerline{\scshape Yan-na Wu}
%\medskip
%{\footnotesize
% %% please put the address of the second author
% \centerline{ Department of  Mathematics, Shanghai Normal University, Shanghai,
%200234, P.R. China}
%   \centerline{Other lines, Springfield, MO 65801-2604, USA}
%} %
%
%\bigskip

%% The name of the associate editor will be entered by an editorial staff
% \centerline{(Communicated by Zhi-Min Zhang
%)}
%%-------------------------------------------------------------------------------------------------------------------------------------------------------------------------------------------------------------------------------------------------------------%%
\begin{abstract}
    In this paper, we propose a novel neural network framework, the Legendre-Kolmogorov-Arnold Network (Legendre-KAN) method, designed to solve fully nonlinear Monge-Amp\`{e}re equations with Dirichlet boundary conditions. The architecture leverages the orthogonality of Legendre polynomials as basis functions, significantly enhancing both convergence speed and solution accuracy compared to traditional methods. Furthermore, the Kolmogorov-Arnold representation theorem provides a strong theoretical foundation for the interpretability and optimization of the network. We demonstrate the effectiveness of the proposed method through numerical examples, involving both smooth and singular solutions in various dimensions. This work not only addresses the challenges of solving high-dimensional and singular Monge-Amp\`{e}re equations but also highlights the potential of neural network-based approaches for complex partial differential equations. Additionally, the method is applied to the optimal transport problem in image mapping, showcasing its practical utility in geometric image transformation. This approach is expected to pave the way for further enhancement of KAN-based applications and numerical solutions of PDEs across a wide range of scientific and engineering fields.%Several novel phenomena are observed for the first time, which open up new avenues for mathematical verification.
\end{abstract}

%%-------------------------------------------------------------------------------------------------------------------------------------------------------------------------------------------------------------------------------------------------------------%%
\section{Introduction}
%方程简介
The Monge-Amp\`{e}re equation is a fully nonlinear elliptic geometric partial differential equation with a broad range of applications, including the classical problem of reconstructing surfaces with prescribed Gaussian curvature, the optimal transport problem with a quadratic cost function, and various physical applications such as reflection, inverse refraction, and others \cite{Benamou1998,Brenier1991,Brix2015,Caffarelli1999,Evans1998,Frisch2002,Karakhanyan2010,Villani2009}. The classical form of the Monge-Amp\`{e}re partial differential equation is given by
\begin{equation}
    \label{eq1}
\begin{aligned}
    \det \left( D^2u \right) =f\left( \emph{\textbf{x}} , u, \nabla u \right),
\end{aligned}
\end{equation}
where $ D^2 u$ denotes the Hessian matrix of an unknown convex function $u$ and $f$ is a given function that depends on $\emph{\textbf{x}} $, $u$ and its gradient $\nabla u$. %$u:{\mathbb{R}}^n \rightarrow \mathbb{R}$, with $n\geq 1$,
%When the domain is $\Omega \subset \mathbb{R} ^2$ \cite{Brenier1991}, the corresponding Hessian determinant can be expressed as:
%\begin{equation}
%    \label{eq:2}
%\begin{aligned}
 %   \det \left( D^2u \right) =\frac{\partial ^2u}{\partial x^2}\frac{\partial ^2u}{\partial y^2}-\left( \frac{\partial ^2u}{\partial x\partial y} \right) ^2
%\end{aligned}
%\end{equation}

In this paper, we propose a Legendre-Kolmogorov-Arnold Network(Legendre-KAN) method to numerically solve the Monge-Amp\`{e}re equation with Dirichlet boundary conditions, which takes the following form:
\begin{equation}
    \label{eq2}
\begin{aligned}
    \left\{ \begin{aligned}
        \det D^2u&=f, \qquad    \mathrm{in} \ \Omega,\\
        u&=g, \qquad    \mathrm{on} \ \partial \Omega,\\
    \end{aligned} \right.
\end{aligned}
\end{equation}
where $\Omega \subset \mathbb{R} ^n$, with $n\geqslant 1$, is an open, bounded and convex set, and
$f$ is a given strictly positive, locally integrable function defined on $\Omega$. The function $g$ represents the Dirichlet boundary condition, defined on $\partial \Omega$. As discussed in \cite{Caffarelli1990,Cristian01,Lions1993}, the solution to \eqref{eq2} is unique only in the cone of convex functions.

There has been increasing interest in the fully nonlinear Monge-Amp\`{e}re equation and its numerical solutions in recent years. However, solving the numerical solution to problem \eqref{eq2} remains challenging due to the fully nonlinear nature of the Monge-Amp\`{e}re operator, which complicates the use of conventional discretization techniques. Furthermore, since convex solutions play a crucial role in many applications, there is a pressing need for effective numerical methods to address this equation. Several numerical methods have been proposed for solving problem \eqref{eq2}, including the wide-stencil finite difference techniques \cite{Froese2011,Oberman2008}, the finite element methods \cite{Awanou2015_1,Awanou2015_2,Brenner2011,Brenner2012,Dean2006,Feng2009_1,Feng2009_2,Feng2013,Gallistl2023,Liu2019,Neilan2013,Neilan2014,Stetter1973}, the spectral methods \cite{Jin2025, Wang2024}, and other collocation methods \cite{lai2023bivariatesplinebasedcollocation}. However, these methods often struggle with high-dimensional problems and singular states. In the emerging field of neural network-based approaches \cite{Nystrom2023}, there are few applications of such methods to problem \eqref{eq2}, particularly in the context of high-dimensional problems and singular states. To address this gap, this paper applies the Legendre-KAN method to obtain numerical solutions of the Monge-Amp\`{e}re equation and to develop a set of practical algorithms for applications in science and engineering.%contributing to the existing literature.

%The Legendre-KAN method demonstrates strong generality in solving optimal transport problems by addressing key limitations of traditional numerical approaches. Unlike finite element methods, which struggle to handle high-dimensional equations efficiently, Legendre-KAN naturally extends to higher dimensions without suffering from the curse of dimensionality. Moreover, while spectral methods provide high accuracy for structured domains, they face significant challenges in handling complex polygonal boundaries, making them unsuitable for optimal transport problems that require flexible domain adaptation. In addition, traditional multilayer perceptrons (MLPs) often require an excessively large number of parameters and layers to achieve comparable accuracy. This not only increases computational cost but also leads to instability in the results, making them less reliable for solving such complex equations. In contrast, Legendre-KAN achieves better convergence with fewer parameters, ensuring both efficiency and stability in high-dimensional and irregular-domain settings.

Recently, a novel neural network architecture, the Kolmogorov-Arnold Network(KAN) has been proposed \cite{Ab2025,gk01,KB2024,liu2024kan,SK2024,YW2025}, demonstrating superior performance in symbolic function fitting compared to traditional Multi-Layer Perceptron(MLP) networks. Based on the Kolmogorov-Arnold theorem, which states that any multivariate function can be represented as a finite combination of univariate functions, KAN is grounded in a solid theoretical framework that enhances its interpretability \cite{Arnold1966,Kolmogorov1961}. This interpretability provides a clear understanding of the network's operational mechanics, aiding in more effective optimization. A key distinction between MLP and KAN is found in their parameter spaces, while the MLP’s parameter space consists solely of weight matrices, KAN’s includes both weight matrices and activation functions. In other words, KAN not only learns the optimal weight matrices but also the optimal activation functions. In this work, we replace the spline functions in KAN with Legendre polynomials, using them as the network's basis functions. As orthogonal polynomials, Legendre polynomials offer high accuracy in function approximation. By combining the strengths of both models, we achieve excellent performance in solving problem \eqref{eq2}.

The primary objective of this paper is to apply the Legendre-KAN method to solve high-dimensional Monge-Amp\`{e}re equations with Dirichlet boundary conditions. Additionally, we consider singular equations to test the generalizability of the proposed method. The main contributions and key features of this study are summarized as follows:

\begin{itemize}
    \item We propose a novel neural network architecture, the Legendre-KAN method, which outperforms traditional MLP networks in solving the Monge-Amp\`{e}re equation.
    \item The method effectively solves the higher-dimensional Monge-Amp\`{e}re equation with Dirichlet boundary conditions, demonstrating its applicability.
    \item We successfully extend the Legendre-KAN method to handle piecewise and weak singularity solutions in both domain and boundary, highlighting its robustness and versatility.
    \item The method is applied to the optimal transport problem in image transportation, showcasing its practical utility.
\end{itemize}

This paper is organized as follows: In Section 2, we introduce the Legendre-KAN architecture and the basic properties of its basis functions. In Section 3, we present appropriate sampling and convergence schemes for the Monge-Amp\`{e}re equation. In Section 4, we provide numerical examples of both smooth and non-smooth solutions in various dimensions to validate the effectiveness of our method. In Section 5, we apply the method to the optimal transport problem. Finally, Section 6 concludes the paper.

%%-------------------------------------------------------------------------------------------------------------------------------------------------------------------------------------------------------------------------------------------------------------%%
\section{Legendre-KAN method}

\subsection{Kolmogorov-Arnold Network}
Unlike MLP, which relies on the universal approximation theorem, KAN is grounded in the Kolmogorov-Arnold representation theorem. This theorem asserts that every continuous multivariate function can be expressed as a finite composition of univariate continuous functions in a two-layer nested structure. The two layers are referred to as the inner and outer functions, respectively. The theorem represents a more constrained yet more general form of Hilbert's thirteenth problem, which is stated as
\begin{equation}
\label{eq:1}
f\left(\emph{\textbf{x}}\right) =f\left( x_1,\cdots,x_n \right) =\sum_{q=0}^{2n}{\Phi _q\left( \sum_{p=1}^n{\phi _{q,p}\left( x_p \right)} \right)},
\end{equation}
where $\Phi _q:\mathbb{R}\rightarrow\mathbb{R}$ are the outer functions and $\phi _{q,p}:[0,1]\rightarrow\mathbb{R}$ are the inner functions. This representation theorem provides a solid theoretical foundation for neural network architecture design, showing that complex multivariate functions can be expressed through compositions of univariate functions. This not only simplifies network's structure design but also facilitates the handling of high-dimensional data. Moreover, the theorem demonstrates that high-dimensional functions can be decomposed into combinations of lower-dimensional ones, which enhances both the model's interpretability and computational efficiency.

Additionally, this theorem can be extended to a general KAN structure with $L$-layer networks of arbitrary width as shown below:
\begin{equation}
\label{eq:2}
y=\mathrm{KAN}\left( \emph{\textbf{x}} \right) =\left( \Phi _{L-1}\circ \Phi _{L-2}\circ \cdots \circ \Phi _0 \right) \left( \emph{\textbf{x}} \right),
\end{equation}
\begin{equation}
\label{eq:3}
\Phi _l=\left( \begin{matrix}
	\phi _{l,1,1}(\cdot )&		\phi _{l,1,2}(\cdot )&		\cdots&		\phi _{l,1,n_l}(\cdot )\\
	\phi _{l,2,1}(\cdot) &		\phi _{l,2,2}(\cdot )&		\cdots&		\phi _{l,2,n_l}(\cdot )\\
	\vdots&		\vdots&		&		\vdots\\
	\phi _{l,n_{l+1},1}(\cdot )&		\phi_{l,n_{l+1},2}(\cdot )&		\cdots&		\phi _{l,n_{l+1},n_l}(\cdot )\\
\end{matrix} \right),
\end{equation}
where $\Phi_l$ represents the function matrix corresponding to the $l^{th}$ KAN layer, $\phi $ denotes the activation functions, $l=0,\cdots,L-1$ is the layer index, and $n_l$ and $n_{l+1}$ represent the number of nodes in the $l$-th and $(l+1)$-th layers, respectively. The KAN architecture is designed to approximate complex functions by combining simpler functions in a hierarchical manner, which is particularly advantageous for high-dimensional problems.

In KAN, the activation function $\phi (x)$  is defined as the sum of a basis function  $b(x)$, similar to residual connections, and a B-spline interpolation function
\begin{equation}
\label{eq:4}
\phi (x)=w_bb(x)+w_s\mathrm{spline}(x).
\end{equation}
In most cases,
\begin{equation}
\label{eq:5}
b(x)=\mathrm{silu}(x)=x/\left( 1+e^{-x} \right),
\end{equation}
the factors $w_b$ and $w_s$ are included primarily to better control the overall magnitude of the activation function. $\mathrm{spline}(x)$ is a cubic B-spline function, defined as
\begin{equation}
\label{eq:6}
\mathrm{spline}(x)=\sum_i{c_i}B_i(x),
\end{equation}
where $B_i(x)$ denotes the $i$-th cubic B-spline basis function, and $c_i$ are the corresponding coefficients. This parameterization ensures  the function is flexible,  capable of approximating a wide range of complex behaviors, while maintaining smoothness.

\subsection{Legendre polynomials}
Legendre polynomials form an orthogonal family of polynomials defined on the interval $[-1,1]$, and  play a crucial role in numerical integration, solving partial differential equations (PDEs), and function approximation. Key properties of these polynomials include orthogonality, recurrence relations, and their connection to Gaussian quadrature formulas. Let $P_n(x)$ denote the Legendre polynomials of degree $n$, and let $\mathcal{P}_N$ represent the set of all Legendre polynomials of degree at most $N$. Below, we present a compilation of essential formulas and properties \cite{CC01,BG,JS02} required in this paper.

First, the Legendre polynomials $P_n(x)$ satisfy the orthogonality relation over the interval $[-1,1]$,
\begin{equation}\label{orthogonality}
  \int_{-1}^{1}P_n(x)P_m(x) dx=\gamma_n \delta_{mn}, \quad  \gamma_n=\frac{2}{2n+1},
\end{equation}
and the three-term recurrence relation
\begin{equation}\label{three-term recurrence relation}
  (n+1)P_{n+1}(x)=(2n+1)x P_{n}(x)-n P_{n-1}(x),\quad n \geqslant 1,
\end{equation}
with initial conditions $P_0(x)=1,~P_1(x)=x$.

Additionally, the Legendre polynomials $P_n(x)$ can be explicitly expressed using Rodrigues' formula
\begin{equation}
\label{eq:9}
P_n(x)=\frac{1}{2^nn!}\frac{d^n}{dx^n}\left( \left( x^2-1 \right) ^n \right).
\end{equation}
This formula provides a direct method for generating the Legendre polynomials $P_n(x)$ and can also be used to derive recursive relationships and explicit forms for the polynomials. These polynomials form a complete basis for the space of polynomials of degree less than or equal to $n$.

The Legendre polynomials, as a type of basis, form a global function approximation space. They are used to model complex data patterns and relationships across various orders, offering more flexibility than B-spline basis functions. In traditional fitting approaches, Legendre polynomials, due to their higher-degree global nature, provide a more consistent and smaller approximation error, especially in regions with discontinuities where B-splines tend to have reduced accuracy. Leveraging these advantages, we propose the Legendre-KAN method for numerically solving high-dimensional Monge-Amp\`{e}re equations with Dirichlet boundary conditions. Experimental results show that the beneficial properties of Legendre polynomials in function fitting are effectively incorporated into KAN, allowing Legendre-KAN to outperform Spline-KAN in terms of fitting accuracy for a variety of complex functions. Furthermore, the polynomial degree can be adjusted, offering greater flexibility in the model's capacity.

\subsection{Architecture of Legendre-KAN Network}
With the increasing attention on the KAN network, several variants have emerged. Researchers have explored replacing the B-spline functions in KAN with orthogonal polynomials to improve its function-fitting performance \cite{SS2024,xu2024fourierkangcf}. In this work, we adopt Legendre polynomials as an alternative to B-spline functions. To ensure that the input functions are mapped to the domain of Legendre polynomials, which is defined on the interval $[-1,1]$, we employ a transformation function $t(x)$ to map the inputs accordingly. This function is defined as
\begin{equation}
\label{eq:10}
x= t(x_{\rm{in}})= \tanh \left( x_{\mathrm{in}} \right),\ x^{*} = t(x^{*}_{in}) = \tanh(x^{*}_{\rm{in}}),
\end{equation}
where $x_{\mathrm{in}}$ represents the input function, and $x^{*}$ is the transfered input function by legendre polynomials. The transformation function $t(x)$ maps the input functions to the interval $[-1,1]$, enabling the effective use of Legendre polynomials for function approximation.

\begin{figure}[!ht]
    \centering
    \includegraphics[width=1\textwidth]{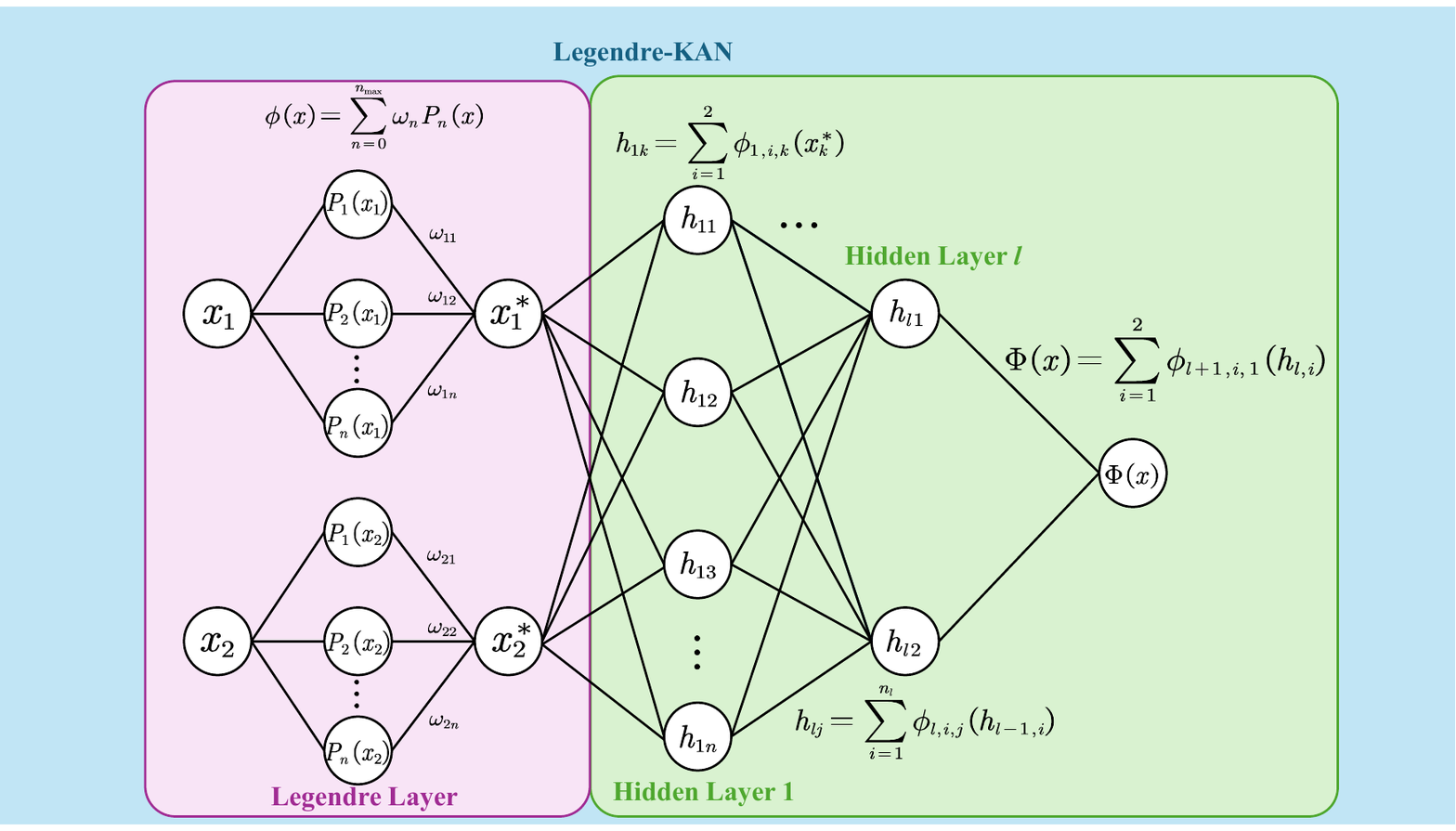} \vspace{-0.7in}
    \caption{Legendre-KAN network}
    \label{fig:1}
\end{figure}

After applying the transformation, the independent variable is fed into the network as
\begin{equation}
\label{eq:11}
\phi \left( x \right) =\sum_{n=0}^{n_{\max}}{\omega _nP_n\left( x \right)},
\end{equation}
where $\omega _n$ represents the weights of the Legendre polynomials. The network architecture is shown in Fig \ref{fig:1}.\vspace{-0.4in}
%%-------------------------------------------------------------------------------------------------------------------------------------------------------------------------------------------------------------------------------------------------------------%%
\section{Optimizer for Monge-Amp\`{e}re equation}% Discretization and solution methods

\subsection{Adaptive moment estimation}
The Adaptive Moment Estimation (Adam) optimizer is a first-order gradient-based optimization algorithm that combines the advantages of the Momentum method and Root Mean Square Propagation (RMSProp) algorithm. It dynamically adjusts the learning rate during training, which accelerates the parameter update process and enhances convergence efficiency.

Adam maintains two moving averages of the gradients, denoted as $m_t$ and $v_t$, which are updated as follows:
\begin{align*}
    m_t &= \beta_1 m_{t-1} + (1 - \beta_1) g_t, \\
    v_t &= \beta_2 v_{t-1} + (1 - \beta_2) g_t^2,
\end{align*}
where $m_t$ and $v_t$ represent the first and second moment estimates, $g_t$ is the gradient of the loss function at time step $t$, $\beta_1$ and $\beta_2$ are hyperparameters controlling the decay rates of the moving averages, with $t$ being the current time step.

To reduce bias in the early stages of training, Adam applies bias correction for $m_t$ and $v_t$,
\begin{equation}
    \label{eq:12}
    \hat{m}_t = \frac{m_t}{1 - \beta_1^t} \quad \text{and} \quad \hat{v}_t = \frac{v_t}{1 - \beta_2^t}.
\end{equation}

Finally, the optimizer updates the parameters $\theta$ using the corrected moment estimates:
\begin{equation}
    \label{eq:13}
    \theta_{t+1} = \theta_t - \eta \frac{\hat{m}_t}{\sqrt{\hat{v}_t} + \epsilon},
\end{equation}
where $\eta$ is the learning rate and $\epsilon$ is a small constant added to prevent division by zero.

\subsection{Modified loss function}

To solve the fully nonlinear Monge-Amp\`{e}re equation with Dirichlet boundary conditions, we decompose the loss function into two components: interior loss and boundary loss, both computed using the mean squared error (MSE) method. For convenience, in Sections 3.2 and 3.3,  we take the domain $\Omega$ as an example in the two-dimensional plane, with the high-dimensional case following a similar approach.

For the interior loss, within the domain $\Omega$, the target equation is
\begin{equation}
    \label{eq:14}
    \det \left( D^2u\left( x,y \right) \right) =f\left( x,y \right), \qquad  \left( x,y \right) \in \Omega,
\end{equation}
where $D^2u\left( x,y \right)$ is the Hessian matrix of $u\left( x,y \right)$, and $f\left( x,y \right)$ is the given function. The interior loss is defined as
\begin{equation}
    \label{eq:15}
    \mathcal{L} _{\mathrm{interior}}=\frac{1}{N_{\mathrm{int}}}\sum_{i=1}^{N_{\mathrm{int}}}{\bigg( \det \left( D^2u_{\theta}\left( x_i,y_i \right) \right) -f\left( x_i,y_i \right) \bigg) ^2},
\end{equation}
where $\left\{ \left( x_i,y_i \right) \right\} _{i=1}^{N_{\mathrm{int}}}$ are the randomly sampled points in the domain $\Omega$, and $N_{\mathrm{int}}$ is the number of interior points.
$u_{\theta}\left( x,y \right)$ is the approximate solution of the Monge-Amp\`{e}re equation obtained by the Legendre-KAN network.

For the boundary loss, on the boundary $\partial \Omega$, the target equation is
\begin{equation}
    \label{eq:16}
    u\left( x,y \right) =g\left( x,y \right), \qquad  \left( x,y \right) \in \partial \Omega.
\end{equation}
The corresponding boundary loss function is
\begin{equation}
    \label{eq:17}
    \mathcal{L} _{\mathrm{boundary}}=\frac{1}{N_{\mathrm{bnd}}}\sum_{j=1}^{N_{\mathrm{bnd}}}{\bigl( u_{\theta}\left( x_j,y_j \right) -g\left( x_j,y_j \right) \bigl) ^2},
\end{equation}
where $\left\{ \left( x_j,y_j \right) \right\} _{j=1}^{N_{\mathrm{bnd}}}$ are the randomly sampled points on the boundary $\partial \Omega$, and $N_{\mathrm{bnd}}$ is the number of boundary points.

The total loss function, combining both the interior loss and boundary loss, is defined as
\begin{equation}
    \label{eq:18}
    \mathcal{L} _{\mathrm{total}}=\lambda \mathcal{L} _{\mathrm{interior}}+ \mathcal{L} _{\mathrm{boundary}},
\end{equation}
where $\lambda$ is a regularization parameter that balances the relative importance of the interior and boundary losses.

\subsection{Adaptive sampling based on residuals}
During  training, an adaptive sampling strategy is employed to dynamically adjust the sampling density, improving the accuracy of the solution, especially in regions with large errors. The algorithm selects sample points exhibiting higher errors, based on the discrepancy between the numerical and analytical solutions, and applies denser sampling in those regions to improve model fitting ability.

Given the numerical solution $u_{\theta}\left( x_i,y_i \right)$ and the analytical solution $u\left( x_i,y_i \right)$, the error at the $i$-th sampling point is defined as
\begin{equation}
    \label{eq:19}
    \mathrm{error}\left( x_i,y_i \right) =\left| u_{\theta}\left( x_i,y_i \right) -u\left( x_i,y_i \right) \right|.
\end{equation}
This error measures the deviation between the numerical and analytical solutions, reflecting the model's accuracy at each sampling point.

In each training cycle, the algorithm computes the error for all sampling points and selects those with larger errors for re-sampling.
Specifically, in the $\left( \mathrm{epoch}+1 \right)$-th iteration, the points with the largest errors are selected for subsequent operations. The index set $\mathcal{I} _{\mathrm{high} \mathrm{error}}$ of  $k$ high-error points is obtained by selecting all indices satisfying
\begin{equation}
    \label{eq:20}
    \mathcal{I} _{\mathrm{high} \mathrm{error}}=\big\{i|~~ \mathrm{error}\left( x_i,y_i \right)\geq \varepsilon \big\},
\end{equation}
where $\varepsilon$ is the error tolerance threshold to be updated. Here, $\mathcal{I} _{\mathrm{high} \mathrm{error}}$ identifies regions requiring additional sampling refinement. %computing the absolute value of the error and i\in\mathrm{num}_{\mathrm{high} \ \mathrm{error} \ \mathrm{samples}}

For the selected high-error points $\left\{ \left( x_i,y_i \right) \right\} _{i\in \mathcal{I} _{\mathrm{high} \mathrm{error}}}$, new sampling points are generated by randomly perturbing their positions.

Let the selected high-error points be denoted as $\left( x_{\mathrm{higherror}}, y_{\mathrm{higherror}} \right)$. The new sampling points, \\$\left( x_{\mathrm{new}},y_{\mathrm{new}} \right)$, are then generated using the following formula
\begin{equation}
    \label{eq:21}
    x_{\mathrm{new}}=x_{\mathrm{higherror}}+\Delta _x,\quad y_{\mathrm{new}}=y_{\mathrm{higherror}}+\Delta _y,
\end{equation}
where $\Delta _x$ and $\Delta _y$ are random perturbations drawn from a uniform distribution $\mathcal{U} \left(-\delta ,\delta \right)$, with $\delta$ controlling the magnitude of the perturbation. Typically, $\delta$ is a small value to ensure that the new sampling points remain close to the high-error points, thus increasing the local density of the training set.

The new points generated through perturbation $\left( x_{\mathrm{new}},y_{\mathrm{new}} \right)$ are added to the existing sampling set $\left\{ \left( x_{\mathrm{interior}},y_{\mathrm{interior}} \right) \right\}$
\begin{align*}
    \mathcal{X} _{\mathrm{interior}}\gets \mathcal{X} _{\mathrm{interior}}\cup \mathcal{X} _{\mathrm{new}},\quad \mathcal{Y} _{\mathrm{interior}}\gets \mathcal{Y} _{\mathrm{interior}}\cup \mathcal{Y} _{\mathrm{new}},
\end{align*}
where $\mathcal{X} _{\mathrm{new}}$ and $\mathcal{Y} _{\mathrm{new}}$ are the new sampling points  generated through perturbation. The approach increases the number of sampling points in regions with larger errors, providing more training information to improve solution accuracy.

After each adaptive sampling step, the number of new sampling points added is controlled by a predefined increment, called the adaptive sample increment. This strategy gradually increases the size of the training dataset, particularly in regions with higher errors. The increment control helps balance training efficiency and computational resource consumption, allowing the network to perform more iterative training in regions where the error is higher.

\begin{figure}[H]
    \centering
    \includegraphics[width=1\textwidth]{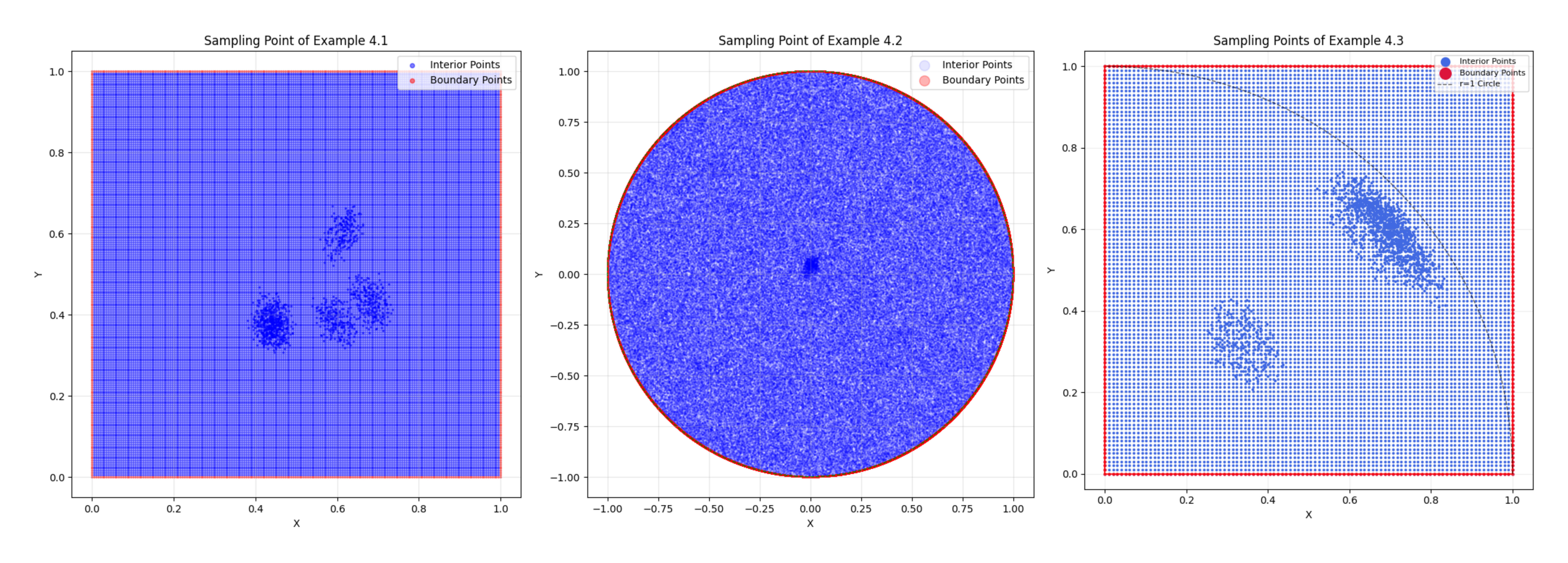} \vspace{-0.4in}
    \caption{Adaptive sampling for the examples}
    \label{fig:adaptive_sampling}
\end{figure}

The core idea of this adaptive sampling strategy is to dynamically adjust the distribution of training points by calculating the residuals and performing local refinement sampling in regions with larger errors. This allows the network to focus on areas with higher error, providing more training information where the solution error is large. This strategy effectively improves the accuracy of the numerical solution and ensures that the neural network focuses on training the regions with the largest errors, thus optimizing the quality of the solution.
Fig\ref{fig:adaptive_sampling} is the examples of the adaptive sampling for the Monge-Ampère equations in the numerical experiments, which can prove that our methods can be used in different types of equations.

%%-------------------------------------------------------------------------------------------------------------------------------------------------------------------------------------------------------------------------------------------------------------%%
\section{ Numerical experiments}
In this section, we present numerical results obtained by solving the Monge-Amp\`{e}re equation in various dimensional settings with the Dirichlet boundary conditions.
Additionally, we compare our results with those produced by MLP applied to the two dimensional non-singular Monge-Amp\`{e}re equation. To evaluate accuracy, we compute two error metrics: the max error, defined as the largest absolute difference between the numerical and exact solutions across all sampling points, and the average error, defined as the mean absolute difference across all points. For clarity, we illustrate the computation using the two-dimensional case, with the same approach extended to higher dimensions. The corresponding formulas are
\begin{align}
    \label{eq:35}
    \rm{error}_{\rm{max}}=\max_{i}{\{ \left| u_{\theta}\left( \emph{x}_i, \emph{y}_i \right) -u\left( \emph{x}_i, \emph{y}_i \right) \right| \}},
    \\ \rm{error}_{\rm{avg}}=\frac{1}{N_{\mathrm{bnd}}+N_{\mathrm{int}}}\sum_{i=1}^{N_{\mathrm{bnd}}+N_{\mathrm{int}}}{\big( \left| u_{\theta}\left( \emph{x}_i, \emph{y}_i \right) -u\left( \emph{x}_i, \emph{y}_i \right) \right| \big)},
\end{align}
where $u_{\theta}$ denotes the numerical solution obtained by the Legendre-KAN network and $u$ represents the exact solution.

\subsection{Example 4.1 Radial Smooth Source in 2D}
We consider the smooth radial function
\begin{equation}
    \label{eq:22}
    u(x,y)=e^{\frac{x^2+y^2}{2}},
\end{equation}
on the unit square $\Omega=(0,1)\times(0,1)$. A straightforward calculation demonstrates that
\begin{equation}
    \label{eq:23}
    f(x,y):= \det (D^2 u(x,y))= (1+x^2+y^2)e^{x^2+y^2}.
\end{equation}
Denote $g$ as the restriction of $u$ to $\partial \Omega$, where $u$ is a convex solution to the Dirichlet problem of \eqref{eq2} with right-hand side $f$ defined by \eqref{eq:23} and boundary values $g$.

\begin{figure}[H]
    \centering
    \includegraphics[width=1\textwidth]{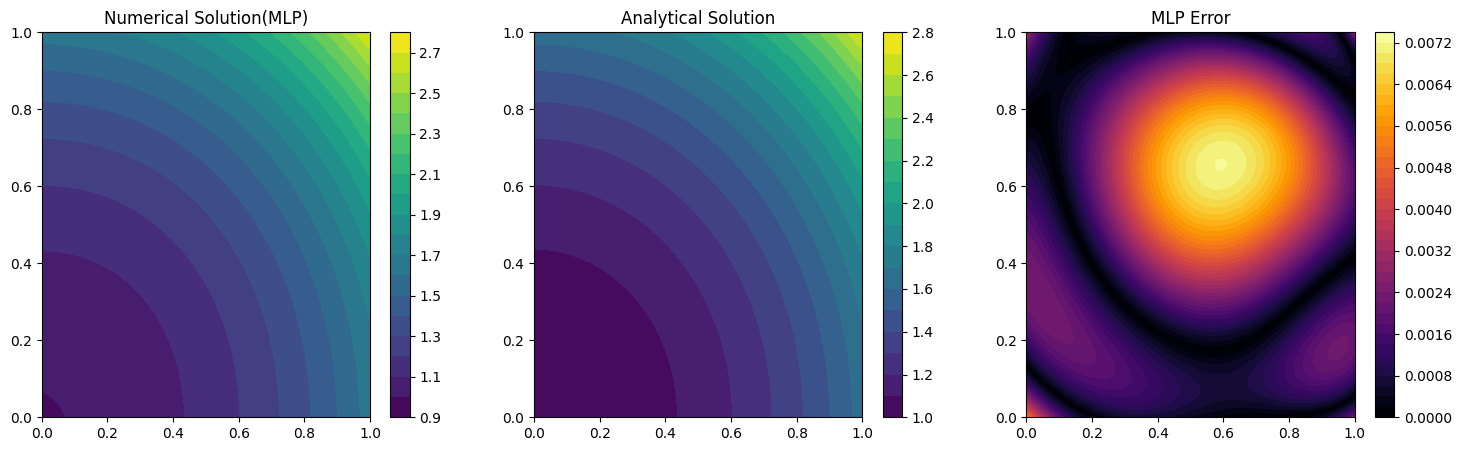}\vspace{-0.2in}
    \caption{ MLP solution for problem \eqref{eq:23}}
    \label{fig:2}
\end{figure}

\begin{figure}[H]
    \centering
    \includegraphics[width=1\textwidth]{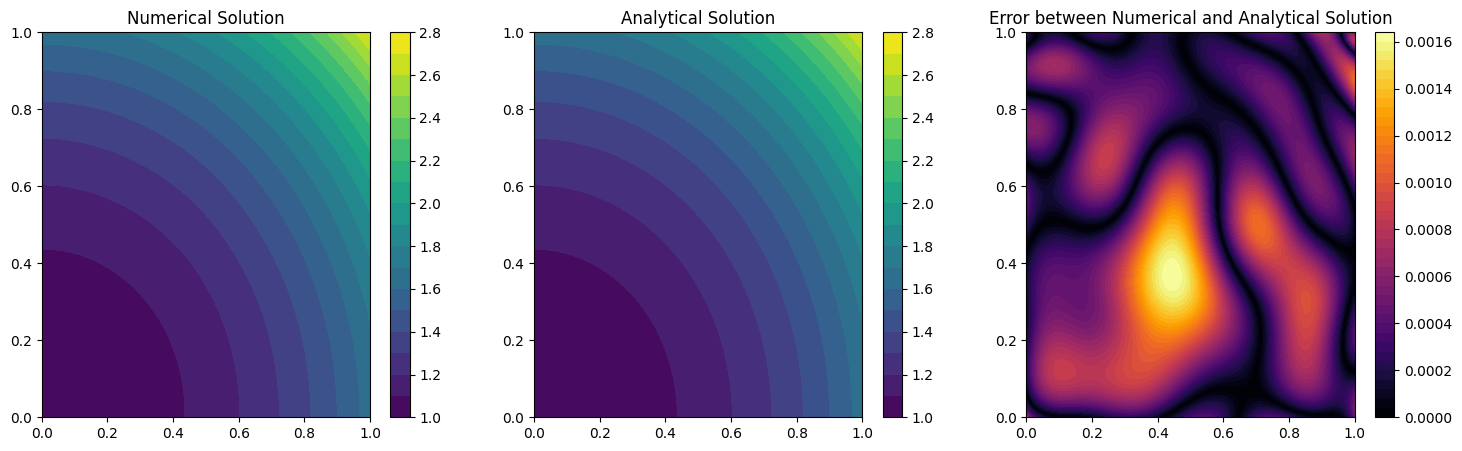}\vspace{-0.2in}
    \caption{Legendre-KAN solution for problem \eqref{eq:23}}
    \label{fig:3}
\end{figure}

\begin{figure}[H]
    \centering
    \includegraphics[width=0.6\textwidth]{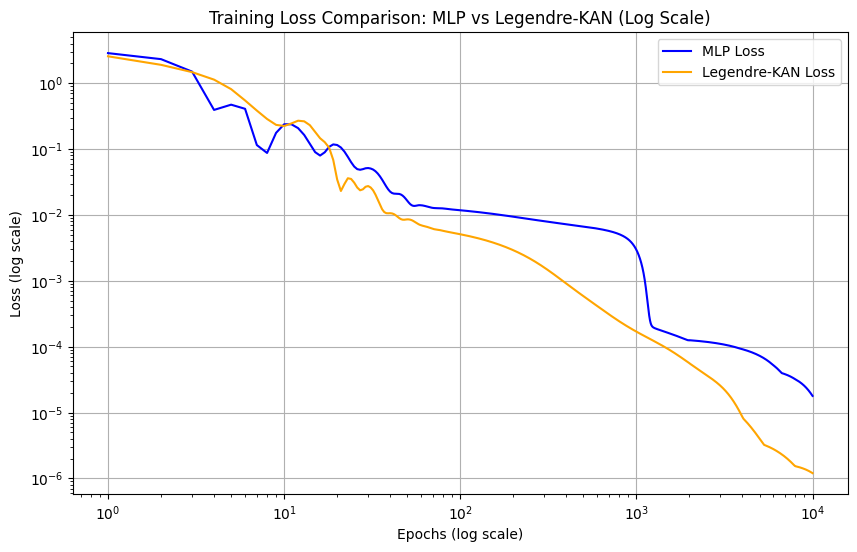}
    \caption{The comparison of the loss function for problem \eqref{eq:23}}
    \label{fig:4}
\end{figure}

In the following, we will deal with inhomogeneous boundary conditions by Legendre-KAN network and MLP network. The parameters of MLP and Legendre-KAN are shown in Table \ref{tab:6}. In Table \ref{tab:1}, we provide a numerical comparison of Legendre-KAN and MLP for problem \eqref{eq:23}. It is evident that our method outperforms MLP in both error values and convergence speed (Fig \ref{fig:2} and \ref{fig:3}) in less parameter counts and simplier network. The loss functions of both methods are illustrated in Fig\ref{fig:4}.

\begin{table}[H]
    \caption{Comparison of MLP and Legendre-KAN Network Architectures}
    \label{tab:6}
    \centering
    \begin{tabular}{lcc}
    \toprule
    \textbf{Attribute} & \textbf{MLP} & \textbf{Legendre-KAN} \\ \midrule
    Number of Hidden Layers & 5 & 2 \\
    Number of Units per Layer & 64 & 8 \\
    Number of Sampling Points & 100 $\times$ 100 & 100 $\times$ 100 \\
    Activation Function & SiLU & Legendre polynomial expansion (degree 6) \\
    Learning Rate & 0.01 & 0.01 \\
    $\lambda$ & 1e-5 & 1e-5 \\
    \bottomrule
    \end{tabular}
\end{table}
\vspace{-0.4in}
\begin{table}[H]
    \caption{Numerical comparison between Legendre-KAN and MLP for problem \eqref{eq:23}}
\label{tab:1}
\centering
\begin{tabular}{lcccc}
    \toprule  % 上边线
    Epoches &  Max Error (MLP) & Average Error (MLP) & Max Error (L-KAN) & Average Error (L-KAN) \\
    \midrule  % 中间线
    2000  & 2.5614e-02  & 7.8799e-03 & 3.0988e-02 & 8.8013e-03\\
    4000  & 2.1003e-02  & 6.2241e-03 & 7.9581e-03 & 2.6861e-03\\
    6000  & 1.4624e-02  & 4.3464e-03 & 2.9458e-03 & 8.5039e-04\\
    8000  & 1.0678e-02  & 3.2886e-03 & 1.6961e-03 & 5.5504e-04\\
    10000 & 7.2107e-03  & 2.3997e-03 & 1.6394e-03 & 4.7582e-04\\
    \bottomrule  % 下边线
    \end{tabular}
\end{table}
\vspace{-0.4in}
\vspace{-0.4in}

Through the above comparison, we observe that our method achieves faster convergence and better convergence performance with a fewer parameters compared to MLP. Building on these results, we aim to further investigate the generality of Legendre-KAN in the context of the Monge-Amp\`{e}re equation. To this end, we will explore different domains and examples involving singular solutions.

\subsection{Example 4.2  Blow-up of the Source Function in different two-dimensional domains}
We consider the function:
\begin{equation}
    \label{eq:24}
    u(x,y)=\left( x^2+y^2\right)^{\alpha}.
\end{equation}
A straightforward calculation demonstrates that
\begin{equation}
    \label{eq:25}
    f(x,y):= \det (D^2 u(x,y))= 4 \alpha^2 (2\alpha -1)(x^2 +y^2)^{2\alpha-1}.
\end{equation}
As $\alpha$ decreases, it's important to note that the singularity of $f$ is getting stronger at the origin. %which lies on the boundary of the domain.

\begin{figure}[H]
    \centering
    \includegraphics[width=0.45\textwidth]{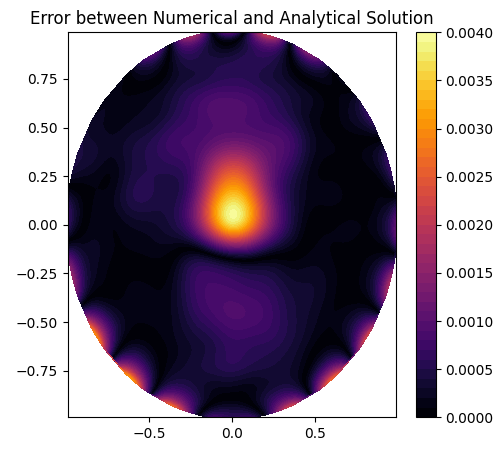}\vspace{-0.2in}
    \caption{Error distribution of the Legendre-KAN solution for problem \eqref{eq:25}}
    \label{fig:5}
\end{figure}

When $\alpha = \frac{5}{3}$ is chosen in problem \eqref{eq:25} with $\Omega = \{ (x,y)|x^2 + y^2=1\}$, the derivative of the function displays a weak singularity at the node $(0, 0)$.
We took 5000 sampling points on the boundary and 400000 sampling points in the area. In Table \ref{tab:2}, we present the error of the numerical solution obtained by the Legendre-KAN method. It is evident that our method can solve the problem effectively (Fig \ref{fig:5}).

\begin{table}[H]
    \caption{Error of the numerical solution for problem \eqref{eq:25} when $\alpha = \frac{5}{3}$}
\label{tab:2}
\centering
\begin{tabular}{lcc}
    \toprule  % 上边线
    Epoches & Max Error & Average Error  \\
    \midrule  % 中间线
    2000  & 1.8029e-02 & 3.4555e-03 \\
    4000  & 9.5146e-03 & 8.3358e-04  \\
    6000  & 6.6757e-03 & 7.9307e-04  \\
    8000  & 4.9622e-03 & 6.4746e-04  \\
    10000 & 3.9409e-03 & 5.2477e-04  \\
    \bottomrule  % 下边线
    \end{tabular}
\end{table}

\begin{figure}[H]
    \centering
    \includegraphics[width=0.45\textwidth]{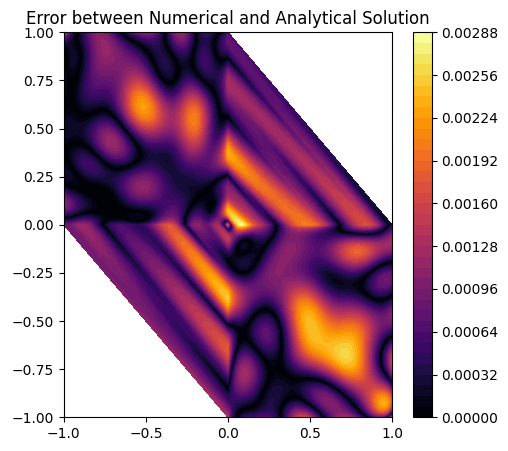}\vspace{-0.2in}
    \caption{Error distribution of the Legendre-KAN solution for problem \eqref{eq:25}}
    \label{fig:6}
\end{figure}

When $\alpha = \frac{3}{4}$ is chosen in problem \eqref{eq:25} with $\Omega=\{(x,y) |-1< x+y <1,~~ \text{and}, -1< x,y <1 \}$, the derivative of the function exhibits stronger singularity at the node $(0, 0)$.
We took 20000 sampling points in the area. In Table \ref{tab:3}, we present the error of the numerical solution obtained by the Legendre-KAN method. The accuracy is also satisfactory (Fig \ref{fig:6}).

\begin{table}[H]
    \caption{Error of the numerical solution for problem \eqref{eq:25} when $\alpha = \frac{3}{4}$}
\label{tab:3}
\centering
\begin{tabular}{lcc}
    \toprule  % 上边线
    Epoches & Max Error & Average Error  \\
    \midrule  % 中间线
    2000  & 9.1863e-03 & 2.1531e-03 \\
    4000  & 5.2631e-03 & 1.2424e-03 \\
    6000  & 4.2815e-03 & 8.0742e-04  \\
    8000  & 3.2548e-03 & 8.2895e-04  \\
    10000 & 2.8084e-03 & 7.3896e-04  \\
    \bottomrule  % 下边线
    \end{tabular}
\end{table}

We conducted numerical experiments on the Monge-Amp\`{e}re equation with singular solutions in different domains. The results demonstrate that the Legendre-KAN method is able to maintain high accuracy in these cases, further highlighting its robustness and generality in handling complex solution structures.

It is worth noting that the Monge-Amp\`{e}re equation has widespread applications in geometric optics and optimal transport theory. For instance, light propagation paths in inhomogeneous media often exhibit piecewise or singularities, and piecewise solutions can effectively model light behavior at medium interfaces or density transition regions. Therefore, studying singular and piecewise solutions is not only theoretically significant but also contributes to advancements in optical design and image reconstruction algorithms.

Moreover, traditional spectral methods often struggle to directly solve problems with irregular boundaries, as they typically rely on structured meshes or specific basis functions. In contrast, neural network-based methods, such as Legendre-KAN, have the distinct advantage of being able to handle equations on arbitrary boundaries by imposing sampling point constraints. This flexibility greatly enhances their applicability to complex geometric domains.

Building on these experimental results, we further investigate examples of the Monge-Amp\`{e}re equation with piecewise solutions. Since piecewise solutions are of great importance in practical problems, we aim to experimentally verify the effectiveness and stability of the Legendre-KAN method in handling such cases.
\subsection{Example 4.3 Piecewise source functions.}
\begin{figure}[H]
    \centering
    \includegraphics[width=0.35\textwidth]{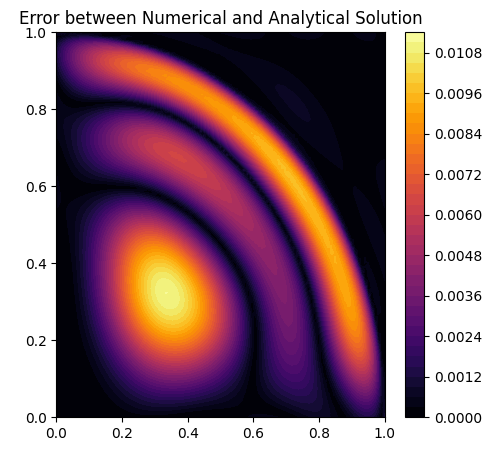}\vspace{-0.2in}
    \caption{Error distribution of the Legendre-KAN solution for problem \eqref{eq:30}}
    \label{fig:8}
\end{figure}

We consider the following problem:
\begin{equation}
    \label{eq:29}
    u\left( r \right) = \begin{cases}
        0,&r\leqslant 1,\\
        r\sqrt{r^2-1}-\ln \left( r+\sqrt{r^2-1} \right) ,\quad &r>1.\\
    \end{cases}
\end{equation}
A straightforward calculation demonstrates that
\begin{equation}
    \label{eq:30}
    f\left( r \right) = \begin{cases}
        0,&		r\leqslant 1,\\
        4,\quad&		r>1,\\
    \end{cases} \
\end{equation}
where $\Omega = (0,1) \times (0,1)$, we took 10000 sampling points in the area.
In Table \ref{tab:7}, we present the error of the numerical solution obtained by the Legendre-KAN method. It is evident that our method can solve the problem effectively (Fig \ref{fig:8}).

\begin{table}[H]
    \caption{Error of the numerical solution for problem \eqref{eq:30}}
\label{tab:7}
\centering
\begin{tabular}{lcc}
    \toprule  % 上边线
    Epoches & Max Error & Average Error  \\
    \midrule  % 中间线
    4000  & 2.2981e-02 & 5.4725e-03 \\
    8000  & 1.7256e-02 & 3.2934e-03 \\
    12000 & 1.3755e-02 & 2.6114e-03  \\
    16000 & 1.1019e-02 & 2.6465e-03  \\
    20000 & 1.1103e-02 & 2.9164e-03  \\
    \bottomrule  % 下边线
    \end{tabular}
\end{table}

\begin{figure}[H]
    \centering
    \includegraphics[width=0.35\textwidth]{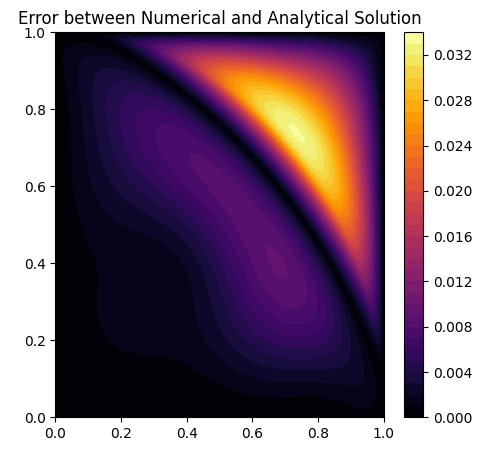}\vspace{-0.2in}
    \caption{Error distribution of the Legendre-KAN solution for problem \eqref{eq:32}}
    \label{fig:9}
\end{figure}

After that, we consider another equation:
\begin{equation}
    \label{eq:31}
    u\left( x,y \right) =\max \left\{ \frac{\left( x-0.5 \right) ^2+\left( y-0.5 \right) ^2}{2},0.08 \right\}.
\end{equation}
A straightforward calculation demonstrates that
\begin{equation}
    \label{eq:32}
    f\left( x,y \right) =\begin{cases}
        1,\quad &\mathrm{if}\ \left( x-0.5 \right) ^2+\left( y-0.5 \right) ^2>0.15^2,\\
        0,&\mathrm{otherwise},\\
    \end{cases}
\end{equation}
where $\Omega = (0,1) \times (0,1)$, we took 10000 sampling points in the area.
In Table \ref{tab:8}, we present the error of the numerical solution obtained by the Legendre-KAN method. The accuracy is also satisfactory (Fig \ref{fig:9}).

\begin{table}[H]
    \caption{Error of the numerical solution for problem \eqref{eq:32}}
\label{tab:8}
\centering
\begin{tabular}{lcc}
    \toprule  % 上边线
    Epoches & Max Error & Average Error  \\
    \midrule  % 中间线
    4000  & 6.4817e-02 & 9.8045e-03 \\
    8000  & 3.8014e-02 & 6.0000e-03 \\
    12000 & 3.4852e-02 & 5.8119e-03  \\
    16000 & 3.3606e-02 & 5.9575e-03  \\
    20000 & 3.3187e-02 & 6.0944e-03  \\
    \bottomrule  % 下边线
    \end{tabular}
\end{table}

In the examples with piecewise solutions, numerical results show that the Legendre-KAN method has surpassed traditional finite element methods in terms of accuracy. This advantage fully demonstrates the potential and superiority of neural networks in handling complex solution structures.
Another remarkable advantage of neural networks is their ability to efficiently solve high-dimensional problems, which is often a major challenge for traditional finite element methods. Thanks to the powerful expressive capability and parameter-sharing characteristics of neural networks, the Legendre-KAN method exhibits good adaptability and stability when extending to higher dimensions.
To further verify the performance of this method in high-dimensional problems, we will present examples of three-dimensional and four-dimensional Monge-Amp\`{e}re equations in the following sections, exploring the accuracy and efficiency of Legendre-KAN in higher dimensions.

\subsection{Example 4.4 Higher-Dimensional Examples.}
We consider the three-dimensional Monge-Amp\`{e}re equation as follows:
\begin{equation}
    \label{eq:26}
    \begin{cases}
        \det \left( D^2u\left( x,y,z \right) \right) =\left( 1+x^2+y^2+z^2 \right) e^{\frac{3\left( x^2+y^2+z^2 \right)}{2}},\quad&		\left( x,y,z \right) \in \Omega,\\
        u\left( x,y,z \right) =e^{\frac{x^2+y^2+z^2}{2}},&		\left( x,y,z \right) \in \partial \Omega, \\
    \end{cases}
\end{equation}
where $\Omega =(0,1) \times (0,1) \times (0,1)$. The exact solution is $u\left( x,y,z \right) =e^{\frac{x^2+y^2+z^2}{2}}$. To solve the problem, we employed the Legendre-KAN method using 125000 sampling points within the domain. The resulting error is presented in Table \ref{tab:4}. As demonstrated by the results, the Legendre-KAN method achieves high accuracy on the 3D Monge-Amp\`{e}re equation.

\begin{table}[H]
    \caption{Error of the numerical solution for problem \eqref{eq:26}}
\label{tab:4}
\centering
\begin{tabular}{lcc}
    \toprule  % 上边线
    Epoches & Max Error & Average Error  \\
    \midrule  % 中间线
    4000  & 3.3850e-02 & 1.1135e-03 \\
    8000  & 1.6960e-02 & 6.4061e-04 \\
    12000 & 1.9277e-02 & 6.1667e-04  \\
    16000 & 1.1706e-02 & 5.1749e-04  \\
    20000 & 9.6383e-03 & 4.8441e-04 \\
    \bottomrule  % 下边线
    \end{tabular}
\end{table}

After that, we consider the four-dimensional Monge-Amp\`{e}re equation as follows:
\begin{equation}
    \label{eq:27}
    \begin{cases}
        \det \left( D^2u\left( x,y,z,w \right) \right) =\left( 1+x^2+y^2+z^2+w^2 \right) e^{\left( x^2+y^2+z^2+w^2 \right)},\quad&		\left( x,y,z,w \right) \in \Omega,\\
        u\left( x,y,z,w \right) =e^{\frac{x^2+y^2+z^2+w^2}{2}},&		\left( x,y,z,w \right) \in \partial \Omega, \\
    \end{cases}
\end{equation}
where $\Omega =(0,1) \times (0,1) \times (0,1)\times (0,1)$. The exact solution is given by $u\left( x,y,z,w \right) =e^{\frac{x^2+y^2+z^2+w^2}{2}}$. We employed the Legendre-KAN method to solve the problem using 390625 sampling points within the domain. The resulting error is reported in Table \ref{tab:5}. These results demonstrate that the Legendre-KAN method also maintains good accuracy when applied to high-dimensional Monge-Amp\`{e}re equation.

\begin{table}[H]
    \caption{Error of the numerical solution for problem \eqref{eq:27}}
\label{tab:5}
\centering
\begin{tabular}{lcc}
    \toprule  % 上边线
    Epoches & Max Error & Average Error  \\
    \midrule  % 中间线
    4000  & 1.3353e-01 & 5.0458e-03 \\
    8000  & 6.7133e-02 & 3.8067e-03 \\
    12000 & 5.7263e-02 & 3.4383e-03 \\
    16000 & 4.4092e-02 & 3.3622e-03  \\
    20000 & 3.8233e-02 & 3.5325e-03 \\
    \bottomrule  % 下边线
    \end{tabular}
\end{table}

These results indicate that the Legendre-KAN method maintains both high accuracy and fast convergence, even in high-dimensional settings.

\subsection{Example 4.5 No Knowing the Exact Solution Case} 
\begin{figure}[H]
    \centering
    \includegraphics[width=0.4\linewidth]{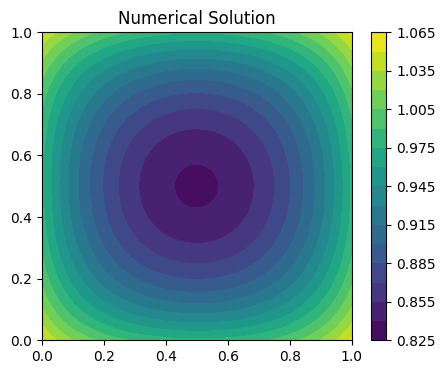}\vspace{-0.2in}
    \caption{The numerical solution of problem \eqref{eq:28} by Legendre-KAN method}
    \label{fig:7}
\end{figure}

We consider the following equation:

\begin{equation}
    \label{eq:28}
    \begin{cases}
        \det \left( D^2u\left( x,y \right) \right) =1,\quad&		\left( x,y \right) \in \Omega,\\
        u\left( x,y \right) =1,&		\left( x,y \right) \in \partial \Omega, \\
    \end{cases}
\end{equation}
with $\Omega =(0,1) \times (0,1)$. We use the Legendre-KAN method to solve the problem. Despite the simplicity of its parameters, problem \eqref{eq:28} lacks smooth classical solutions due to the presence of corners in $\Omega$, despite it being the unit square.
In Fig \ref{fig:7}, we show the numerical solution of problem \eqref{eq:28} by Legendre-KAN method. Numerical results affirm the effectiveness of our method, which doesn't require knowing the exact solution to the Monge-Amp\`{e}re equation.

Through this example, we can see that although we cannot provide the exact solution to the Monge-Amp\`{e}re equation, the Legendre-KAN method can still effectively solve it. This not only demonstrates the robustness and generality of the method but also lays a solid foundation for solving optimal transport problems. With this approach, we can map images onto any plane, enabling flexible transportation and transformation.

\section{Applications of the Optimal Transport Problem}
In optimal transport theory, constructing the optimal transport map is a fundamental problem. To address this challenge, we leverage the celebrated Brenier theorem \ref{thm:1}, which provides a powerful theoretical foundation for the uniqueness and construction of optimal transport maps \cite{Brenier1991}. The theorem establishes that, under certain conditions, the optimal transport map can be expressed as the gradient of a convex potential function. This result is of great significance in both the mathematical theory of optimal transport and its practical applications, especially in fields such as computer vision and image reconstruction, where it significantly improves algorithmic efficiency and accuracy.

Before formally stating the theorem, let us briefly discuss the problem setting: Given two probability measures \(\mu\) and \(\nu\) on \(\mathbb{R}^d\), where \(\mu\) is absolutely continuous with respect to the Lebesgue measure and \(\nu\) has finite second moments, our goal is to find an optimal transport map that minimizes the transportation cost, typically given by the squared Euclidean distance. Brenier's theorem precisely characterizes the properties of this map, as stated below:

\begin{theorem}{Brenier Theorem (1991)}
\label{thm:1}
    Let \(\mu\) and \(\nu\) be probability measures on \(\mathbb{R}^d\), where \(\mu\) is absolutely continuous with respect to the Lebesgue measure and \(\nu\) has finite second moments. Then there exists a unique (up to \(\mu\)-a.e. equivalence) convex function \(\varphi: \mathbb{R}^d \to \mathbb{R}\) such that:
    \begin{itemize}
        \item The gradient map \(\nabla \varphi\) pushes \(\mu\) forward to \(\nu\), i.e.,
        \[
        \nabla \varphi \# \mu = \nu.
        \]
        \item This map \(\nabla \varphi\) is the unique optimal transport map between \(\mu\) and \(\nu\), minimizing the quadratic cost:
        \[
        \int_{\mathbb{R}^d} |\textbf{x} - \nabla \varphi(\textbf{x})|^2 \, d\mu(\textbf{x}).
        \]
    \end{itemize}
\end{theorem}

In our numerical experiments, we have demonstrated that the Legendre-KAN method exhibits high accuracy and fast convergence when solving the Monge-Amp\`{e}re equation. To further illustrate its practical value in optimal transport applications, we consider a special case in this section. Based on Brenier's theorem, which asserts that the optimal transport map for convex cost functions is given by the gradient of a convex potential, we select a smooth and convex function as the right-hand side of the equation. This choice not only facilitates performance validation but also enables direct application to common 2D image processing tasks. Supported by this theoretical foundation, our method demonstrates robust computational potential and broad applicability in computer vision, image reconstruction, and other engineering domains.

We consider the Monge-Amp\`{e}re equation arising in the optimal transport problem:
\begin{equation}
    \label{eq:33}
    \begin{cases}
        \det \left( D^2u\left( \emph{\textbf{x}} \right) \right) =f\left( \emph{\textbf{x}} \right) ,\quad&		\emph{\textbf{x}}\in V,\\
        \nabla u\left( V \right) =V,&		u\,\,\mathrm{is} \ \mathrm{convex} \ \mathrm{over} \ V.\\
    \end{cases}
\end{equation}
Given an input image $I:V \subset \mathbb{R} ^3\rightarrow \mathbb{R} ^3$, where $V=(0,1) \times (0,1)$ represents the 2D domain of the image, and the image has RGB channels.
The goal is to construct a transformation map to warp the image coordinates while satisfying certain constraints.
The right-hand $f(\emph{\textbf{x}})$ is a given source function that defines the desired geometric property of the transformation.

\subsection{Example 5.1} We consider the following Monge-Amp\`{e}re equation with an optimal transport boundary condition:

\begin{equation}
\label{eq:ot1}
    \begin{cases}
        \det \left( D^2u\left( \emph{\textbf{x}}  \right) \right) = 1,\quad&		\emph{\textbf{x}} \in V,\\
        \nabla u\left( V \right) =V,&		u\,\,\mathrm{is} \ \mathrm{convex} \ \mathrm{over} \ V.\\
    \end{cases}
\end{equation}
To solve this equation, we apply the Legendre-KAN method and obtain the optimal transport map $u(\emph{\textbf{x}})$, which can be used to warp the input image $I$ to achieve the desired geometric transformation.
\begin{figure}[H]
    \centering
    \includegraphics[width=0.75\textwidth]{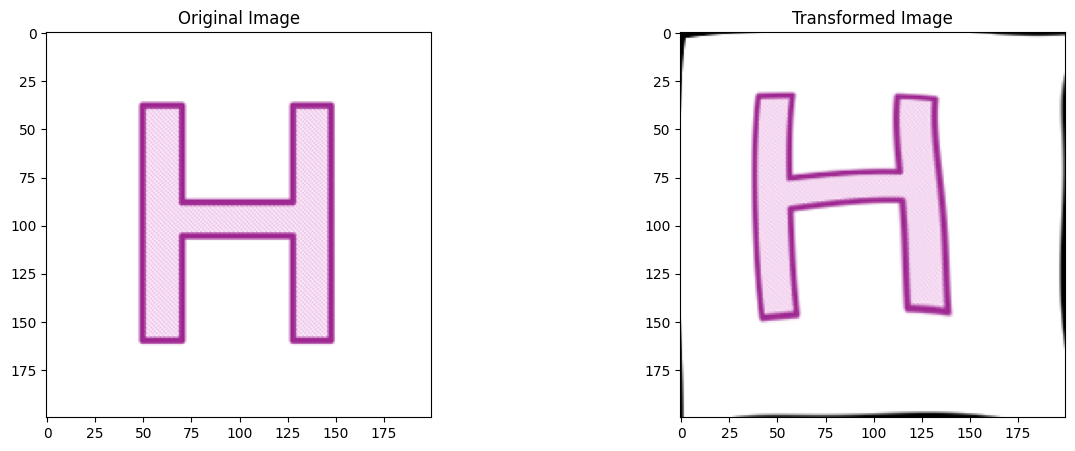}\vspace{-0.2in}
    \caption{The transformation of the input image by equation \ref{eq:ot1}}
    \label{fig:ot1}
\end{figure}\vspace{-0.1in}
Having computed a simple example of the optimal transport boundary for the Monge-Amp\`{e}re equation, we now turn our attention to a more illustrative case. By modifying the objective function, we aim to demonstrate that the Legendre-KAN method can be employed to map an image onto a convex surface defined by the target function. In the following discussion, we present an example where a fisheye image is optimally transported onto a nearly flat plane, showcasing the method’s capability to achieve geometrically meaningful image transformations.

\subsection{Example 5.2} To address the fisheye distortion, we consider the following Monge-Amp\`{e}re equation with optimal transport boundary condition:
\begin{equation}
\label{eq:ot2}
\begin{cases}
    \det \left( D^2u\left( \emph{\textbf{x}}  \right) \right) = f\left( \emph{\textbf{x}}  \right) / \rm{A}(\emph{f}~),\quad&		\emph{\textbf{x}} \in V,\\
    \nabla u\left( V \right) =V,&		u\,\,\mathrm{is} \ \mathrm{convex} \ \mathrm{over} \ V,\\
\end{cases}
\end{equation}
where $f\left( \emph{\textbf{x}}  \right)$ is a given function that defines the desired geometric property of the transformation, and $\rm{A}(\emph{f}~)$ denotes the whole image element. To correct the fisheye distortion and transform the image into a more natural perspective, we choose the following Gaussian function:
\begin{equation}
    f(\emph{\textbf{x}} )=A \exp \left(-\frac{\left(x-x_0\right)^2+\left(y-y_0\right)^2}{2 \sigma^2}\right),
\end{equation}
where $(x_0,y_0)$ is the center of the fisheye image, $\sigma$ is the standard deviation of the Gaussian function. The Legendre-KAN method is used to solve this equation and obtain the optimal transport map $u(x,y)$, which is subsequently used to warp the input image $I$ and achieve the desired geometric correction.

\begin{figure}[H]
    \centering
    \includegraphics[width=0.75\textwidth]{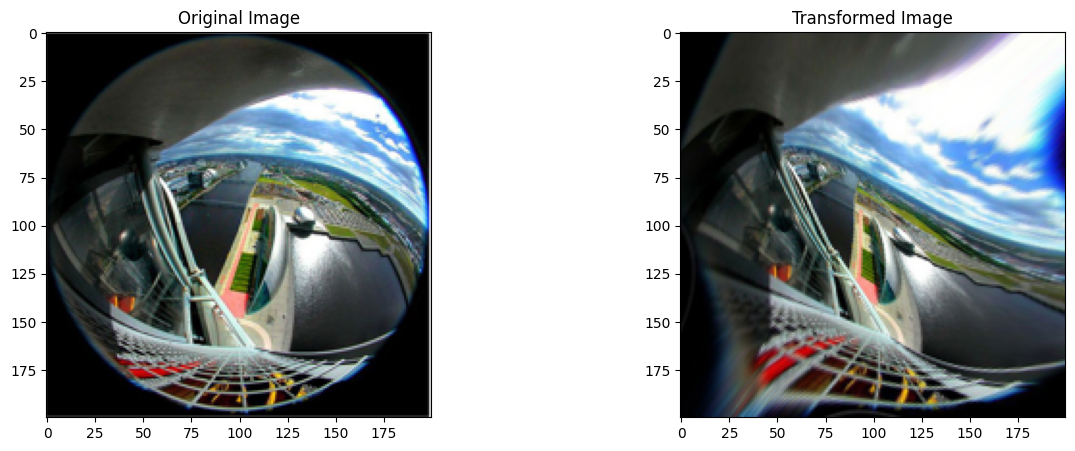}\vspace{-0.2in}
    \caption{The transformation of the input image by equation \ref{eq:ot2}}
    \label{fig:ot2}
\end{figure}
After the fisheye image transformation, we demonstrate that the Legendre-KAN method exhibits strong generality in addressing optimal transport problems in image processing. This versatility allows it to handle a variety of image transformation tasks, such as mapping complex image structures onto different geometries, including curved or non-Euclidean spaces.

Furthermore, this method can be extended to similar applications in fields like image registration, texture mapping, and even shape analysis, where optimal transport plays a crucial role. These applications benefit from the flexibility of Legendre-KAN in dealing with arbitrary boundary conditions and high-dimensional problems, making it a powerful tool for a wide range of computational imaging challenges.
%%-------------------------------------------------------------------------------------------------------------------------------------------------------------------------------------------------------------------------------------------------------------%%
 %\newpage
 \section{Concluding Remarks}
 In this work, we propose the Legendre-KAN method to solve the Monge-Amp\`{e}re equation under Dirichlet boundary conditions, demonstrating its effectiveness in both smooth and non-smooth settings. By leveraging the Kolmogorov-Arnold representation theorem, our method utilizes Legendre polynomials as basis functions, providing a theoretically sound and numerically robust framework. The key contributions of this study include:
\begin{enumerate}
    \item Developing a Legendre-KAN architecture that improves convergence speed and accuracy compared to traditional MLP networks.
    \item Designing an adaptive sampling strategy that focuses on regions with large errors, dynamically refining the training set to enhance solution accuracy.
    \item Demonstrating the versatility of Legendre-KAN in solving high-dimensional and singular Monge-Amp\`{e}re equations through various numerical experiments.
    \item The method is applied to the optimal transport problem in image mapping, demonstrating its practical effectiveness in geometric image transformation.
\end{enumerate}
Our results show that the Legendre-KAN consistently outperforms traditional approaches such as MLP in terms of both error reduction and computational efficiency. The proposed method proves particularly effective for high-dimensional Monge-Amp\`{e}re equations, maintaining accuracy and stability while overcoming challenges associated with nonlinearity and convexity. These findings lay a promising foundation for future theoretical analysis and formal proofs of the method's properties. Moreover, as our approach is straightforward and easy to implement, we expect the results presented in this paper to inspire further applications within the scientific and engineering communities.
\\

%\vskip 0.2in
\noindent {\bf Data Availability} The datasets generated during and/or analysed during the current study are available from the corresponding author on reasonable request.
\vskip 0.1in
\noindent {\bf Declarations}\\
\noindent {\bf Conflict of interest} The authors declared that they have no conflict of interest.
%%-------------------------------------------------------------------------------------------------------------------------------------------------------------------------------------------------------------------------------------------------------------%%
%\newpage

\bibliographystyle{plain}
\end{document}